\documentclass[reqno]{amsart}
\usepackage[margin = 1.5in]{geometry}
\usepackage{amsmath, amssymb, amsthm, fancyhdr, verbatim, graphicx, mathtools, commath}
\usepackage{enumitem}
\usepackage[all]{xy}
\usepackage{mathrsfs}
\usepackage{hyperref}
\usepackage{amsbsy}
\usepackage{booktabs}
\usepackage{float}
\usepackage{bm}
\usepackage{xurl}

\numberwithin{equation}{section}
\newcommand{\seqnum}[1]{\href{https://oeis.org/#1}{\rm \underline{#1}}}

\newcommand{\nc}{\newcommand}
\nc{\renc}{\renewcommand}

\nc{\bZ}{\mathbb{Z}}
\nc{\bN}{\mathbb{N}}

\nc{\cS}{\mathcal{S}}
\nc{\cA}{\mathcal{A}}
\nc{\cB}{\mathcal{B}}
\nc{\cH}{\mathcal{H}}

\nc{\flr}[1]{\lfloor #1 \rfloor}
\nc{\ceil}[1]{\lceil #1 \rceil}

\newtheorem{theorem}{Theorem}[section]
\newtheorem{lemma}[theorem]{Lemma}
\newtheorem{proposition}[theorem]{Proposition}
\newtheorem{corollary}[theorem]{Corollary}

\theoremstyle{remark}
\newtheorem{remark}[theorem]{Remark} 
\newtheorem{definition}[theorem]{Definition}
\newtheorem{example}[theorem]{Example}

\makeatletter
\def\th@remark{%
  \thm@headfont{\bfseries}%
  \normalfont 
  \thm@preskip \thm@preskip 
  \thm@postskip\thm@preskip
}

\title[Counting Second-Order Schreier Sets with Minimum 2]{Counting Second-Order Schreier Sets with Minimum 2}

\author{Dongho Kim}
\address{Institute of Mathematics, Academia Sinica, Taipei, Taiwan}
\email{dkim@as.edu.tw}

\begin{document}

\begin{abstract}
    A finite nonempty subset of the positive integers is first-order Schreier if the minimum element of that set is greater than or equal to the size of the set. A nonempty subset of the positive integers $X$ is second-order Schreier if $X$ is a union of $d$ first-order Schreier sets for some $d \leq \min X$. We study a combinatorial relationship between the count of second-order Schreier sets and the Fibonacci sequence. Specifically, we show that the second-order Schreier set with minimum $2$ and maximum $n+4$ is one more than the $n^{\text{th}}$ term of the convolution of the Fibonacci numbers $(F_1, F_2, \ldots)$ with $(F_2, F_3, \ldots)$ (given by Sequence A023610 in the OEIS).
\end{abstract}

\maketitle

\tableofcontents
\section{Introduction}
Banach and Saks originally asked whether every weakly convergent sequence in a Banach space admits a subsequence whose Ces\`aro means converge in norm~\cite{BanSaks}. Schreier~\cite{SchreierEinGZ} provided a counterexample to this question by introducing Schreier sets, which are defined as follows.
\begin{definition}\label{def: Schreier set}
    A finite nonempty set $X \subset \bZ_{>0}$ is a \textit{first-order Schreier set} if
    \begin{equation}\label{eq: def of first-order Schreier}
        \abs{X} \ \leq \ \min X.
    \end{equation}
    For convenience, we say a first-order Schreier set is \textit{Schreier}, and we denote the collection of Schreier sets by $\cS_1$ and refer to it as the \textit{first Schreier family}. We also include the empty set, so
    \begin{equation}
        \cS_1 \ \coloneq \ \{ \varnothing\} \ \cup \ \bigl\{ X \subset \bZ_{>0} : \abs{X} \leq \min X \bigr\}.
    \end{equation}
    We denote the collection of Schreier sets with $n$ as their largest element by $\cS_1(n)$.
\end{definition}
Schreier sets and their properties also appear in combinatorics, specifically in Ramsey theory~\cite{far2004}. Accordingly, Schreier sets have notable properties and applications in both functional analysis and combinatorics, and one direction of interest is the generalization of Schreier sets to higher orders.
\begin{definition}\label{def: 2nd order Schreier, nth order Schreier}
    A finite nonempty set $X \subset \bZ_{>0}$ is a \textit{second-order Schreier set} if
    \begin{equation}
        X \ = \ \bigcup_{i=1}^d X_i
    \end{equation}
    such that each $X_i \in \cS_1$ and $d \leq \min X$. We denote the collection of second-order Schreier sets by $\cS_2$ and refer to it as the \textit{second Schreier family}, so
    \begin{equation}
        \cS_2 \ \coloneq \ \{ \varnothing \} \cup \Biggl\{ \bigcup_{i=1}^d X_i: X_i \in \cS_1 \ \text{and} \ \min X \geq d \Biggr\}.
    \end{equation}
    We denote the collection of second-order Schreier sets with $a$ as their minimum element and $b$ as their maximum element by
    \begin{equation}
        \cS_2(a,b) \ \coloneq \ \{X \in \cS_2 : \min X = a, \ \max X = b \}.
    \end{equation}
    More generally, a nonempty set $X \subset \bZ_{>0}$ is an \textit{$n^{\text{th}}$-order Schreier set} if
    \begin{equation}
        X \ = \ \bigcup_{i=1}^d X_i
    \end{equation}
    such that each $X_i \in \cS_{n-1}$ and $1 \leq d \leq \min X$. We denote the collection of $n^{\text{th}}$-order Schreier sets by $\cS_n$ and refer to it as the \textit{$n^{\text{th}}$ Schreier family}. Note that the empty set is an element of $\cS_n$ for all $n \in \bZ_{>0}$.
\end{definition}
Another direction to explore is the count of Schreier sets. Bird~\cite{out-of-the-norm} previously showed that the number of Schreier sets with $n$ as their largest element is equal to the $n^{\text{th}}$ Fibonacci number.
\begin{theorem}[\cite{out-of-the-norm}]\label{thm: counting Schreier sets} 
    Let ${(F_n)}_{n \geq 0}$ be the Fibonacci sequence with $F_0=0$, $F_1=1$, and
    \[
    F_{n+2} \ = \ F_{n+1} + F_n \quad \text{for } n \geq 0.
    \]
    Then, the number of Schreier sets with $n$ as their largest element is $F_n$, i.e.,
    \begin{equation}
        \abs{\cS_1(n)} \ = \ F_n.
    \end{equation}
\end{theorem}
Several authors have proved Theorem~\ref{thm: counting Schreier sets} using formula-based and bijective  arguments~\cite{out-of-the-norm,Chu2019,Chu2021,Chu2026arxiv,ChuGeng,ChuIrmark,ChuMillerXiang,ChuVasseur}.

We focus on the combinatorial properties of Schreier sets, and in particular, we study the count of $\cS_2(2,n)$, the second-order Schreier sets whose minimum is $2$ and maximum is $n$. There are two reasons why we choose to study $\cS_2(2,n)$. First, the second-order Schreier sets are the immediate generalization of the first-order Schreier sets and those with minimum $2$ are the simplest to count. Second, because an arbitrary set in $\cS_2(2,n)$ is a union of two sets (one of which is possibly empty) in $\cS_1$, we presume that the count of $\cS_2(2,n)$ can be expressed in terms of the Fibonacci numbers. The following table presents the values of $\abs{\cS_2(2,n)}$ for $n = 2, \ldots, 12$.
\begin{table}[H]
    \centering
    \setlength{\tabcolsep}{12pt}
    \begin{tabular}{cc} 
        \toprule
        $n$ & $ \abs{\cS_2(2,n)}$ \\
        \midrule
        2 & 1 \\
        3 & 1 \\
        4 & 2 \\
        5 & 4 \\
        6 & 8 \\
        7 & 16 \\ 
        8 & 31 \\
        9 & 59 \\
        10 & 110 \\
        11 & 202 \\
        12 & 366 \\
        \bottomrule
    \end{tabular}
    \caption{Count of (min $2$, max $n$) second-order Schreier sets.}\label{tab: values of S_2(2,n)}
\end{table}
While the count of $\cS_2(2,n)$ seems unremarkable, we observe an interesting pattern when we consider the following sequence.
\begin{definition}\label{def: convolution of Fibonacci numbers}
    Let ${(T_n)}_{n\geq 0}$ be the sequence given by \textit{the convolution of ($F_1, F_2, \ldots$) and ($F_2, F_3, \ldots$)}, where $F_n$ is the $n^{\text{th}}$ Fibonacci number. That is,
    \begin{equation}
        T_n \ = \ \sum_{i=0}^{n} F_{i+1} F_{n+2-i}.
    \end{equation}
    While the sequence begins at $T_0=1$, we refer to $T_n$ as the $n^{\text{th}}$ term of the convolution of $(F_1, F_2, \ldots)$ and $(F_2, F_3, \ldots)$.
\end{definition}
Note that $T_n$ appears in OEIS~\cite{oeis} as sequence \seqnum{A023610}~\cite{A023610}, and the following table gives the values of $T_n$.
\begin{table}[H]
    \centering
    \setlength{\tabcolsep}{16pt}
    \begin{tabular}{cc} 
        \toprule
        $n$ & $ T_n$ \\
        \midrule
        0 & 1 \\
        1 & 3 \\
        2 & 7 \\
        3 & 15 \\
        4 & 30 \\
        5 & 58 \\ 
        6 & 109 \\
        7 & 201 \\
        8 & 365 \\
        \bottomrule
    \end{tabular}
    \caption{Convolution of $(F_1, F_2, \ldots)$ and $(F_2, F_3, \ldots)$.}\label{tab: values of T_n}
\end{table}
The values in Table~\ref{tab: values of S_2(2,n)} and Table~\ref{tab: values of T_n} point to a deeper relationship between $\abs{\cS_2(2,n)}$ and $T_n$. Formalizing this relationship is our main result.
\begin{theorem}\label{thm: count of min 2 second-order Schreier sets}
    For $n \geq 0$, the number of second-order Schreier sets with minimum $2$ and maximum $n+4$ is one more than $T_n$. That is, 
    \begin{equation}
        \abs{\cS_2(2,n+4)} \ = \ 1 + T_n.
    \end{equation}
\end{theorem}
To build up to this statement, we present an equivalent characterization of second-order Schreier sets in Section~\ref{sec: preliminaries}, which we use to split an arbitrary second-order Schreier set into a disjoint union of Schreier sets $A_1, \ldots, A_r$ ordered \textit{successively}, i.e., $A_1 < \cdots <A_r$ with $A_i < A_j$ if and only if $\max A_i < \min A_j$. We also present a recurrence relation for $T_n$, which is key in the proof of Theorem~\ref{thm: count of min 2 second-order Schreier sets}. Following this, we provide an alternative proof of Theorem~\ref{thm: counting Schreier sets} in Section~\ref{sec: counting first Schreier}. 
The proof technique for Theorem~\ref{thm: counting Schreier sets} is used to count a collection of Schreier sets given in the following definition.
\begin{definition}\label{def: cH(a,b) and H(a,b)}
    Let $a, b \in \bZ_{>0}$ with $a \leq b$. We define the collection of Schreier sets whose minimum is at least $a$ and maximum is exactly $b$ as
    \begin{equation}
        \cH(a,b) \ \coloneq \ \bigl\{X \in \cS_1(b): \min X \geq a \bigr\}.
    \end{equation}
    We define $H$ as the count of such sets, i.e.,
    \begin{equation}
        H(a,b) \ \coloneq \ \abs{\cH(a,b)} \ = \ \# \bigl\{X \in \cS_1(b): \min X \geq a \bigr\}.
    \end{equation}
\end{definition}
Based on the definition of $\cH(a,b)$, we observe the following recurrence relation for $H(a,b)$ given $a,b \in \bZ_{>0}$ with $a \leq b$.
\begin{equation}
    H(a+1, b+2) \ = \ H(a+1, b+1) + H(a,b).
\end{equation}
In Section~\ref{sec: counting min 2 second Schreier}, we show that $\cS_2(2,n)$ can be expressed in terms of $\cH(a,n)$ for varying $a$, which helps us count $\cS_2(2,n)$ with $H(a,n)$. We then present a second recurrence relation for $H(a,b)$, this time involving a series of $H(a,b)$ terms, and we relate this recurrence relation to the recurrence relation for $T_n$, which we use to prove Theorem~\ref{thm: count of min 2 second-order Schreier sets}.

\section{Preliminaries}\label{sec: preliminaries}
We introduce two tools that we use in the proof of Theorem~\ref{thm: count of min 2 second-order Schreier sets}. The first is an equivalent characterization of second-order Schreier sets, and the second is a recurrence relation for $T_n$, the convolution of Fibonacci numbers $(F_1, F_2, \ldots)$ and $(F_2, F_3, \ldots)$.
\vspace{1em}
\subsection{Equivalent characterization of second-order Schreier sets}~\label{subsec: alternative definition of 2nd order Schreier}
\newline
Recall from Definition~\ref{def: 2nd order Schreier, nth order Schreier} that a nonempty set $X \subset \bZ_{>0}$ is second-order Schreier if $X$ can be written as
\[
X \ = \ \bigcup_{i=1}^d X_i,
\]
where $X_i \in \cS_1$ for all $i$ and $d \leq \min X$. Equivalently, we can define second-order Schreier sets as a disjoint union of first-order Schreier sets that are ordered in a certain way.
\begin{lemma}\label{lem: alternative definition of 2nd order Schreier}
    A nonempty set $X \subset \bZ_{>0}$ is a second-order Schreier set if and only if
    \begin{equation}
        X \ = \ \bigsqcup_{i=1}^{r} A_i,
    \end{equation}
    where $A_i \in \cS_1$ is nonempty for all $i$, $A_i < A_{i+1}$ for all $i$, and $r \leq \min A_1$. Note that $A_i < A_{i+1}$ means $\max A_i < \min A_{i+1}$.
\end{lemma}
\begin{proof}
    We note that $A_i < A_{i+1}$ for all $i$ means the union of the $A_i$'s is disjoint. Also, $\min A_1 = \min X$, so the condition $r \leq \min A_1$ is equivalent to $r \leq \min X$.
    \begin{enumerate}[label=(\roman*)]
        \item We first prove the forward direction. Given nonempty $X \in \cS_2$, we can write
        \begin{equation}
            X \ = \ \bigcup_{i=1}^d X_i,
        \end{equation}
        with $X_i \in \cS_1$ and $d \leq \min X$. Let $Y_1 \coloneq X_1$, and for each $i >1$, define
        \begin{equation}
            Y_i \ \coloneq \ X_i \setminus \bigcup_{k < i} X_k.
        \end{equation}
        By definition, the $Y_i$'s are pairwise disjoint and their union is $X$. We also note each nonempty $Y_i$ is Schreier because $\abs{Y_i} \leq \abs{X_i} \leq \min X_i \leq \min Y_i$. We now discard any empty $Y_i$'s and reindex the remaining $Y_i$'s so that $\min Y_1 < \cdots < \min Y_r$ for some $r \leq d$ (note $r$ may equal $d$).

        Let us now define $N_0 \coloneq 0$ and
        \begin{equation}
            N_i \ \coloneq \ \sum_{k=1}^i \abs{Y_k} \quad \text{for } i = 1, \ldots, r.
        \end{equation}
        With this, we index the elements of $X$ in ascending order:
        \begin{equation}\label{eq: indexing elements of X}
            X \ = \ \{x_1 < x_2 < \cdots < x_{N_r} \},
        \end{equation}
        and define
        \begin{equation}\label{eq: defining A_i}
            A_i \ \coloneq \ \{x_{N_{i-1}+1}, \ldots, x_{N_i} \} \quad \text{for each } i.
        \end{equation}
        From~\eqref{eq: indexing elements of X} and~\eqref{eq: defining A_i}, we have
        \[
            X \ = \ \bigsqcup_{i=1}^r A_i,
        \]
        and $A_i < A_{i+1}$ for all $i$. It remains to show that each $A_i$ is Schreier. For $i=1$, we have
        \[
            \min A_1 \ = \ \min Y_1 \ \geq \ \abs{Y_1} \ = \ \abs{A_1}.
        \]
        For $i > 1$, we note that
        \begin{equation}
            \# \{x \in X : x < \min Y_i \} \ \leq \ \sum_{k<i} \abs{Y_k} \ = \ N_{i-1}
        \end{equation}
        by the ordering of the $Y_i$'s. Then, we know from~\eqref{eq: defining A_i} that
        \begin{equation}
            \abs{A_i} \ = \ \abs{Y_i} \ \leq \ \min Y_i \ \leq \ x_{N_{i-1} + 1} \ = \ \min A_i,
        \end{equation}
        so $A_i$ is Schreier. Because $r \leq d \leq \min X = \min A_1$, we are done.

        \item The reverse direction is trivial by taking each $X_i$ to be $A_i$.
    \end{enumerate}
\end{proof}

Recall from Definition~\ref{def: 2nd order Schreier, nth order Schreier} that we defined $\cS_2(a,b)$ to be the collection of second-order Schreier sets with minimum element $a$ and maximum element $b$. We observe from examples that the characterization of second-order Schreier sets in Lemma~\ref{lem: alternative definition of 2nd order Schreier} is consistent with Definition~\ref{def: 2nd order Schreier, nth order Schreier}.
\begin{example}
    Consider, for example, $\cS_2(2, 5)$, the collection of second-order Schreier sets with minimum $2$ and maximum $5$. We observe that
    \[
    \cS_2(2,5) \ = \ \bigl\{ \{2,5\}, \{2,3,5\}, \{2,4,5\}, \{2,3,4,5\} \bigr\}.
    \]
    It is clear that each of these sets can be written as $A_1 \sqcup A_2$ such that $A_1, A_2 \in \cS_1$ and $A_1 < A_2$. For example,
    \[
    \{2,3,4,5\} \ = \ \{2,3\} \sqcup \{4,5\}.
    \]
    For $\cS_2(2,8)$, we observe that
    \[
    \cS_2(2,8) \ = \ \{ A \subsetneq [2,8]: \min A= 2, \ \max A = 8 \}.
    \]
    All strict subsets of $[2,8]$ containing both $2$ and $8$ are second-order Schreier sets with minimum $2$ and maximum $8$. However, $[2,8]$ is not in $\cS_2(2,8)$ since the set cannot be written as the union of two Schreier sets. Using Lemma~\ref{lem: alternative definition of 2nd order Schreier}, we again note there exists no way to write $[2,8]$ as a disjoint union of two Schreier sets.
\end{example}
\begin{remark}
    While Lemma~\ref{lem: alternative definition of 2nd order Schreier} is restricted to second-order Schreier sets, this characterization naturally extends to $n^{\text{th}}$-order Schreier sets for $n > 2$, allowing them to be expressed as a disjoint union of $(n-1)^{\text{th}}$-order Schreier sets ordered in the same way.
\end{remark}
\vspace{1em}
\subsection{Recurrence relation for $T_n$}~\label{subsec: recurrence relation for T_n}
\newline
We present a recurrence relation for $T_n$, the convolution of $(F_1, F_2, \ldots)$ and $(F_2, F_3, \ldots)$. Recall from Definition~\ref{def: convolution of Fibonacci numbers} that $T_n$ is given by the formula
\[
T_n \ = \ \sum_{i=0}^n F_{i+1} F_{n+2-i}.
\]
The Fibonacci relation $F_{n+2} = F_{n+1} + F_n$ is used to derive a recurrence relation for $T_n$.
\begin{proposition}[\cite{A023610}]\label{prop: recurrence relation for T_n}
    For $n \geq 0$, we have the following recurrence relation
    \begin{equation}\label{eq: recurrence relation for T_n}
        T_{n+2} \ = \ F_{n+4} + T_{n+1} + T_n.
    \end{equation}
\end{proposition}
We present a proof here for completeness.
\begin{proof}
    Observe that
    \begin{align}
        F_{n+4} + T_{n+1} + T_{n} \
        &= \ \bigl(F_{n+3} + F_{n+2} \bigr) + \sum_{i=0}^{n+1} \bigl(F_{i+1} F_{n+3-i} \bigr) + \sum_{j=0}^n \bigl(F_{j+1} F_{n+2-j} \bigr) \nonumber \\[2pt]
        &= \ F_{n+3} F_{2} + F_{n+2} F_{1} + F_{n+2} F_2 + \sum_{i=0}^n F_{i+1} \bigl(F_{n+3-i} + F_{n+2-i} \bigr) \nonumber \\[2pt]
        &= \ F_{n+3} F_2 + F_{n+2} F_3 + \sum_{i=0}^n F_{i+1} F_{n+4-i} \nonumber \\
        &= \ \sum_{i=0}^{n+2} F_{i+1} F_{n+4-i} \nonumber \\[7pt]
        {} &= \ T_{n+2}.
    \end{align}
\end{proof}

\section{First-order Schreier sets}\label{sec: counting first Schreier}
Before counting the second-order Schreier sets that have $2$ as their minimum element, we present an alternative method to count the number of (first-order) Schreier sets. This method helps us count the number of Schreier sets whose minimum element is at least $a$ and maximum element is $b$, which, in turn, is useful in counting the second-order Schreier sets.\vspace{1em}
\subsection{An alternative proof}~\label{subsec: alternative proof}
\newline
We recall that a nonempty set $X \subset \bZ_{>0}$ is Schreier if it satisfies~\eqref{eq: def of first-order Schreier}, i.e.,
\[
\abs{X} \ \leq \ \min X.
\]
\begin{example}
    Let us consider some examples and non-examples of Schreier sets. Consider the following sets
    \[
    A \ = \ \{5,6,7\} \quad \text{and} \quad B \ = \ [10,12] \cup [14,21]
    \]
    As $\abs{A} = 3 \leq 5 = \min A$, the set $A$ is Schreier. However, we note that $B$ is not Schreier, as $\abs{B} = 11 \not \leq 10 = \min B$.
\end{example}
Previous work done by Bird~\cite{out-of-the-norm} presented several methods to count the Schreier sets that have $n$ as their largest element, denoted by $\cS_1(n)$. To give an alternative method, we first determine the maximum size of a Schreier set.
\begin{lemma}\label{lem: Schreier set max size}
    The maximum size of a Schreier set with \(n\) as its largest element is \( \ceil{n/2} \) or, equivalently, \(\flr{(n+1)/2} \). Furthermore, there exists a Schreier set of size $\flr{(n+1)/2}$.
\end{lemma}
\begin{proof}
    Let $X \subset \bZ_{>0}$ be a Schreier set with $n$ as its largest element. Then, $\abs{X} \leq \min X \leq n$. We want to show that $\abs{X} \leq \flr{(n+1)/2}$.

    Suppose $\abs{X} = k$. Then, $\min X \geq k$, which means that $X \subset [k, n]$, so $\abs{X} \leq n-k+1$. From this, we have
    \begin{equation}
        k \ \leq \ n-k+1, 
    \end{equation}
    or, equivalently, $k \leq (n+1)/2$. Since $k$ must be an integer, we conclude that
    \begin{equation}
        k \ \leq \ \biggl\lfloor \frac{n+1}{2} \biggr\rfloor.
    \end{equation}
    It remains to show there exists a Schreier set of size $\flr{(n+1)/2}$ for all $n$. We consider two cases depending on the parity of $n$.
    \begin{enumerate}[label=(\roman*)]
        \item Let \(n\) be even. Consider the set
        \[
        A \ \coloneq \ \{n/2, \ldots, n\} \setminus \{n-1 \}.
        \]
        There are $(n/2) + 1$ elements in the set \( \{n/2, \ldots, n\} \), so removing $\{n-1\}$ gives us a set of size $n/2$. We observe
        \[
        \abs{A} \ = \ \frac{n}{2} \ \leq \ \min A,
        \] 
        so $A$ is Schreier.

        \item Now, let $n$ be odd. Then, \( \flr{(n+1)/2} = (n+1)/2 \in \bZ \). There exists a Schreier set of size $(n+1)/2$ given by
        \[
        A' \ \coloneq \ \biggl\{ \frac{n+1}{2}, \frac{n+3}{2}, \dots, \frac{n+n}{2} = n \biggr\}.
        \]
        Note that $\abs{A'} = (n+1)/2 = \min A'$, so $A'$ is Schreier.
    \end{enumerate}
\end{proof}
With Lemma~\ref{lem: Schreier set max size}, we are ready to count $\cS_1(n)$, the collection of Schreier sets with $n$ as their largest element. We know by Theorem~\ref{thm: counting Schreier sets} that
\[
\abs{\cS_1(n)} \ = \ F_n.
\]
The alternative proof relies on individually counting the elements of $\cS_1(n)$ of size $k$ for $1 \leq k \leq \flr{(n+1)/2}$. Note the range of $k$ comes from the Schreier set $\{n \}$ and Lemma~\ref{lem: Schreier set max size}. After this work was completed, the author learned that the same proof appears in~\cite{Chu2026arxiv}.

\begin{proof}[Proof of Theorem~\ref{thm: counting Schreier sets}]
    We note that
    \begin{equation}\label{eq: alternative counting Schreier 1}
        \abs{\cS_1(n)} \ = \ \sum_{k=1}^{\flr{(n+1)/2}} \# \{X \in \cS_1(n) : \ \abs{X} = k \}.
    \end{equation}
    By definition, the minimum element of an arbitrary Schreier set of size $k$ is greater than or equal to $k$. Thus, if $X \in \cS_1(n)$ is of size $k$, then $X = A \cup \{n\}$, where $A$ is a $(k-1)$-subset of $[k, n-1]$. As such, 
    \begin{equation}\label{eq: number of max n, size k Schreier}
        \# \{ X \in \cS_1(n) : \ \abs{X} = k \} \ = \ \binom{(n-1) - k +1}{k-1} \ = \ \binom{n-k}{k-1},
    \end{equation}
    which means
    \begin{equation}\label{eq: alternative counting Schreier 2}
        \abs{\cS_1(n)} \ = \ \sum_{k=1}^{\flr{(n+1)/2}} \binom{n-k}{k-1} \ = \sum_{m=0}^{\flr{(n-1)/2}} \binom{n-m-1}{m} \ = \ F_n,
    \end{equation}
    where the last equality is a well-known combinatorial identity for the Fibonacci numbers.
\end{proof}
\begin{example}
    We verify our proof for $n=6$. Observe that
    \[
    F_{6} \ = \ \sum_{m=0}^{\flr{(n-1)/2}} \binom{n-m-1}{m} \ = \ \binom{5}{0} + \binom{4}{1} + \binom{3}{2},
    \]
    and the elements of $\cS_1(6)$ are
    \[
    \bigl\{ \{6\} \bigr\} \cup \bigl\{ \{2,6\}, \{3,6\}, \{4, 6\}, \{5, 6\} \bigr\} \cup \bigl\{ \{3, 4, 6\}, \{3, 5, 6\}, \{4, 5, 6 \} \bigr\}.
    \]
    It is clear there is one element of size $1$, four elements of size $2$, and three elements of size $3$, which corresponds exactly to the binomials of $F_6$. 
\end{example}
\vspace{1em}
\subsection{Counting a subset of Schreier sets}~\label{subsec: counting a subset of Schreier sets}
\newline
This alternative method of counting $\abs{\cS_1(n)}$ given in~\eqref{eq: alternative counting Schreier 1} and~\eqref{eq: alternative counting Schreier 2} helps us count the number of Schreier sets with a minimum element of at least $a$ and a fixed maximum element of $b$. Recall from Definition~\ref{def: cH(a,b) and H(a,b)} that
\[
\cH(a,b) \ = \ \bigl\{X \in \cS_1(b): \min X \geq a \bigr\},
\]
the collection of Schreier sets with minimum at least $a$ and maximum exactly $b$, and
\[
H(a,b) \ = \ \abs{\cH(a,b)},
\]
the count of such sets. We now present a formula for $H(a,b)$.
\begin{lemma}\label{lem: counting min geq a, max b Schreier}
    Let $a, b \in \bZ_{>0}$ with $a \leq b$. For $a \geq 2$, we have
    \begin{equation}\label{eq: counting H for a geq 2}
        H(a,b) \ = \ \sum_{j=0}^{\min (b-a, a-2)} \binom{b-a}{j} + \sum_{k=0}^{\flr{(b-2a+1)/2}} \binom{b-a-k}{a-1+k}.
    \end{equation}
\end{lemma}
\begin{proof}
    We observe that the Schreier sets with minimum element at least $a$ and maximum element exactly $b$ is given by the following disjoint union
    \begin{equation}\label{eq: union of min geq a, max b schreier}
        \{X \in \cS_1(b): \min X \geq a > \abs{X} \} \ \sqcup \ \{X \in \cS_1(b): \min X \geq a \ \text{and} \ \abs{X} \geq a\}.
    \end{equation}
    As such, to count $H(a,b)$, it suffices to count the two sets of~\eqref{eq: union of min geq a, max b schreier}.
    
    The first set in~\eqref{eq: union of min geq a, max b schreier} counts the elements of $\cS_1(b)$ whose size is less than $a$ and whose minimum element is at least $a$. Thus, if $X$ is in the first set, then
    \begin{equation}
        X \ = \  Y \cup \{b\},
    \end{equation}
    where $Y$ is a $j$-subset of $[a,b-1]$ such that
    \begin{equation*}
        j \ \leq \ \abs{[a, b-1]} \quad \text{and} \quad j \ \leq \ a-2.
    \end{equation*}
    The condition $j \leq a-2$ is a consequence of $a > \abs{X}$, and as we set $a \geq 2$, there is no issue of $j$ being negative. From this, we have
    \begin{align}
        \# \{X \in \cS_1(b): \min X \geq a > \abs{X} \} \
        &= \ \sum_{j=0}^{\min (b-a, a-2)} \binom{(b-1)-a+1}{j} \nonumber\\
        &= \ \sum_{j=0}^{\min (b-a, a-2)} \binom{b-a}{j} \label{eq: counting min geq a, max b schreier, size less than a}.
    \end{align}
    We now count the second set of~\eqref{eq: union of min geq a, max b schreier}. It follows by definition that all Schreier sets of size at least $a$ must have minimum element at least $a$. Therefore,
    \begin{equation}
        \{X \in \cS_1(b): \min X \geq a \ \text{and} \ \abs{X} \geq a\} \ = \ \{X \in \cS_1(b): \ \abs{X} \geq a\}.
    \end{equation}
    From~\eqref{eq: number of max n, size k Schreier}, we have
    \begin{equation}\label{eq: counting min geq a, max b schreier, size at least a}
        \# \{X \in \cS_1(b): \ \abs{X} \geq a \} \ = \ \sum_{\ell=a}^{\flr{(b+1)/2}} \binom{b-\ell}{\ell-1} \ = \ \sum_{k=0}^{\flr{(b-2a+1)/2}} \binom{b-a-k}{a-1+k},
    \end{equation}
    where $k \leq \flr{(b-2a+1)/2}$ is equivalent to the condition that $b-a-k \geq a-1+k$. Combining~\eqref{eq: counting min geq a, max b schreier, size less than a} and~\eqref{eq: counting min geq a, max b schreier, size at least a} gives~\eqref{eq: counting H for a geq 2}.
\end{proof}
\begin{remark}
    We defined $H(a,b)$ as the number of Schreier sets with minimum at least $a$ and maximum exactly $b$. This may seem unnatural, as one may expect to instead count the number of Schreier sets with minimum exactly $a$ and maximum exactly $b$. However, we see in Section~\ref{sec: counting min 2 second Schreier} that $H(a,b)$ is useful in counting the second-order Schreier sets.
\end{remark}
One interesting consequence of Lemma~\ref{lem: counting min geq a, max b Schreier} arises when the largest element is less than twice the minimum element. In this case, we have a simple identity for $H(a,b)$.
\begin{corollary}\label{cor: counting min geq a, max b schreier for b<2a}
    For $b < 2a$, we have the equality
    \begin{equation}
        H(a,b) \ = \ 2^{b-a}.
    \end{equation}
\end{corollary}
\begin{proof}
    Recall by definition that $1 \leq a \leq b$. When $a=1$, the value of $b$ necessarily equals $1$ to satisfy $b < 2a$. Since $H(1,1) = 1 = 2^{1-1}$, it remains to consider when $a \geq 2$. We can then solve for $H(a,b)$ using~\eqref{eq: counting H for a geq 2}.

    We separate our proof into two cases.
    \begin{enumerate}
        \item Suppose $b-a < a-1$. Then, $\flr{(b-2a+1)/2} < 0$ and $\min (b-a, a-2) = b-a$, so
        \begin{equation}
        H(a,b) \ = \ \sum_{j=0}^{b-a} \binom{b-a}{j} \ = \ 2^{b-a}.
        \end{equation}

        \item Suppose $b-a = a-1$. Then, $\flr{(b-2a+1)/2} = 0$ and $\min (b-a, a-2) = a-2$. As such,
        \begin{align}
            H(a,b) \ 
            &= \ \sum_{j=0}^{a-2} \binom{b-a}{j} + \sum_{k=0}^{0} \binom{b-a-k}{a-1+k} \nonumber \\[3pt]
            & = \ \biggl( \binom{b-a}{0} + \cdots + \binom{b-a}{b-a-1} \biggr) + \binom{b-a}{b-a} \nonumber \\[3pt]
            & = \ 2^{b-a}. \label{eq: substituting b-a=a-1}
        \end{align}
    \end{enumerate}
\end{proof}
From the definition of $H(a,b)$, we also have the following corollary.
\begin{corollary}
    For $b \geq 2$, we have the equality
    \begin{equation}
        \abs{\cS_1(b)} = H(1,b) = H(2,b).
    \end{equation}
\end{corollary}
\begin{proof}
    The equality $\abs{\cS_1(b)} = H(1,b)$ is clear by definition. As $\{1\}$ is the only Schreier set that contains $1$, we have the equality $H(1,b) = H(2,b)$.
\end{proof}
Though we do not use the following proposition, we still present the formula for counting instead the number of Schreier sets with minimum exactly $a$ and maximum exactly $b$.
\begin{proposition}
    For $2 \leq a < b$,
    \begin{equation}
        \# \{X \in \cS_1(b): \min X = a \} \ = \ \sum_{j \geq 0}^{\min (b-a-1, a-2)} \binom{b-a-1}{j}.
    \end{equation}
\end{proposition}
\begin{proof}
    With the assumption that $a \geq 2$, note that an arbitrary set in
    \[
    \{X \in \cS_1(b): \min X = a \}
    \]
    can be written as
    \[
    \{a\} \cup Y \cup \{b\}
    \]
    where $Y$ is a $j$-subset of $[a+1, b-1]$ such that $j \leq \abs{[a+1, b-1]}$ and $j \leq a-2$. From this, we have
    \[
    \# \{X \in \cS_1(b): \min X = a \} \ = \ \sum_{j\geq 0} \binom{((b-1) - (a+1) + 1)}{j} \ = \ \sum_{j \geq 0} \binom{b-a-1}{j},
    \]
    where $j \leq \min (b-a-1, a-2)$.
\end{proof}
\vspace{1em}
\subsection{A recurrence relation for $H(a,b)$}~\label{subsec: recurrence relation for H}
\newline
We use the definition of $\cH(a,b)$ to establish a recurrence relation for $H(a,b)$, which is useful in counting the sets in $\cS_2(2, n)$.
\begin{theorem}\label{thm: recurrence relation for H}
    Let $a, b \in \bZ_{>0}$ with $a \leq b$. Then,
    \begin{equation}\label{eq: recurrence relation for H(a,b)}
        H(a+1, b+2) \ = \ H(a+1, b+1) + H(a, b).
    \end{equation}
\end{theorem}
\begin{proof}
    Recall from Definition~\ref{def: cH(a,b) and H(a,b)} and Lemma~\ref{lem: Schreier set max size} that
    \begin{equation}
    H(a+1,b+2) \ = \ \# \{X \in \cS_1(b+2) : \ \min X \geq a+1 \} \ = \ \abs{\cH(a+1,b+2)}
    \end{equation}
    and observe $\cH(a+1, b+2)$ is given by the disjoint union
    \begin{equation}\label{eq: disjoint union for cH(a+1,b+2)}
        \renewcommand{\arraystretch}{1.2}
        \begin{array}{lccc}
            {} & \{X \in \cH(a+1,b+2): b+1 \notin X\} & \sqcup & \{X \in \cH(a+1,b+2): b+1 \in X\} \\
            \eqcolon & \cA_1 & \sqcup & \cA_2
        \end{array}
    \end{equation}
    Therefore, to show~\eqref{eq: recurrence relation for H(a,b)} holds, it suffices to count the two sets of~\eqref{eq: disjoint union for cH(a+1,b+2)}.
    \begin{enumerate}[label=(\roman*)]
        \item There exists a bijection between $\cA_1$ and $\cH(a+1, b+1)$. Let 
        \[
        \varphi: \cA_1 \ \longrightarrow \ \cH(a+1, b+1)
        \]
        be defined by replacing $b+2$ with $b+1$ in a given set $X \in \cA_1$. Then, $\varphi(X)$ is well-defined, as $b+2 \notin \varphi(X)$, $b+1 \in \varphi(X)$, and
        \begin{equation}
            \min \varphi(X) \ = \ \min \bigl\{ b+1, \min X \bigr\} \ \geq \ a+1.
        \end{equation}
        We know $\varphi$ is injective, as $b+1 \notin X$. For surjectivity, observe that if $Y \in \cH(a+1, b+1)$, then
        \begin{equation}
            Y \ = \ \varphi \bigl((Y \cup \{b+2\}) \setminus \{b+1\} \bigr).
        \end{equation}
        We know
        \begin{equation}
            (Y \cup \{b+2 \}) \setminus \{b+1\} \ \in \ \cA_1
        \end{equation}
        because $b+1$ is not in the set, the minimum is at least that of $Y$, and the maximum is $b+2$. As such, 
        \begin{equation}
            \# \{X \in \cH(a+1, b+2): b+1 \notin X \} \ = \ H(a+1, b+1).
        \end{equation}

        \item There exists a bijection between $\cA_2$ and $\cH(a,b)$. Let the map
        \[
        \psi: \cH(a,b) \ \longrightarrow \ \cA_2
        \]
        be defined by
        \[
        \psi(X) \ = \ \{x+1: x \in X \} \cup \{b+2\},
        \]
        which means both $b+1$ and $b+2$ are in $\psi(X)$. Observe that
        \begin{equation}
            \min\psi(X) \ = \ 1 + \min X \ \geq \ a+1,
        \end{equation}
        and recalling the definition of Schreier sets, we note that $\psi(X)$ can contain $1$ more element than $X$. As such, the union of $\{x+1: x \in X\}$ with $\{b+2\}$ is still Schreier, so $\psi(X)$ is well-defined in $\cA_2$. We know $\psi$ is injective, as $b+1 \notin X$. For surjectivity, observe that if $Y \in \cA_2$, then
        \begin{equation}
            Y \ = \ \psi( \{y - 1 : y \in Y \} \setminus \{b+1 \} ).
        \end{equation}
        We know
        \begin{equation}
            \{y - 1 : y \in Y \} \setminus \{b+1 \} \ \in \ \cH(a,b)
        \end{equation}
        because the minimum is at least $(a+1)-1=a$ and removing $b+1$ ensures the set is Schreier with maximum $b+1-1=b$. As such,
        \begin{equation}
            \# \{X \in \cH(a+1,b+2): b+1 \in X\} \ = \ H(a,b).
        \end{equation}
    \end{enumerate}
\end{proof}

\section{Second-order Schreier sets with minimum 2}\label{sec: counting min 2 second Schreier}
We see from~\eqref{eq: alternative counting Schreier 2} that the count of $\cS_1(n)$ is equal to $F_n$, where $F_n$ is the Fibonacci sequence. We study the count of $\cS_2(2,n+4)$, the collection of second-order Schreier sets with minimum element $2$ and maximum element $n+4$, which we show is $1+T_n$, or one more than the convolution of $(F_1, F_2, \ldots)$ and $(F_2, F_3, \ldots)$.
\vspace{1em}
\subsection{Counting $\cS_2(2,n)$ with $H$}~\label{subsec: counting cS_2(2,n) with H}
\newline
We first introduce a way to count $\cS_2(2,n)$ using the function $H(a,b)$, which counts the number of Schreier sets with minimum element at least $a$ and maximum element $b$.
\begin{lemma}\label{lem: count of min 2 second-order Schreier using H}
    Let $n \geq 4$. Then, we have
    \begin{equation}
        \abs{\cS_2(2,n)} \ = \ 1 + \sum_{k=4}^{n} H(k, n)
    \end{equation}
\end{lemma}
\begin{proof}
    For $n \geq 4$, every $X$ in $\cS_2(2,n)$ contains $2$ and $n$, so $\abs{X} \geq 2$. As such, we can separately count $\abs{\cS_2(2,n)}$ based on their cardinality as follows:
    \begin{equation}\label{eq: first disjoint union of cS_2(2,n)}
        \renewcommand{\arraystretch}{1.2}
        \begin{array}{llccc}
            \cS_2(2,n) & = & \{X \in \cS_2(2,n): \abs{X} = 2\} & \sqcup & \{X \in \cS_2(2,n): \abs{X} \geq 3\} \\
            {} & \eqcolon & \cB_1 & \sqcup & \cB_2
        \end{array}
    \end{equation}
    The only set in $\cB_1$ is $\{2,n\}$, so $\abs{\cB_1} = 1$. To count $\cB_2$, we partition it by $\min (X \setminus \{2\})$, which is at least $3$ and at most $n-1$.
    \begin{equation}
        \cB_2 \ = \ \bigsqcup_{k=3}^{n-1} \cB_2^{k}, \quad \text{where} \quad \cB_2^{k} \coloneq \{X \in \cB_2 : \min(X \setminus \{2\}) = k\}.
    \end{equation}
    For fixed $k \in \{3, \ldots, n-1\}$, there exists a bijection between $\cB_2^k$ and $\cH(k+1, n)$. Consider the map
    \[
    \pi_k: \cB_2^k \ \longrightarrow \ \cH(k+1, n)
    \]
    defined by
    \begin{equation}
        \pi_k(X) \ = \ X \setminus \{2, k\}.
    \end{equation}
    Because $\abs{X} \geq 3$ and $\min X =2$, we know from Lemma~\ref{lem: alternative definition of 2nd order Schreier} that $X$ is the disjoint union of two Schreier sets ordered successively:
    \[
    X \ = \ A_1 \sqcup A_2,
    \]
    where $2 \in A_1$ and $n \in A_2$. We now study two cases.
    \begin{enumerate}[label=(\roman*)]
        \item If $k \in A_1$, then $A_1$ necessarily equals $\{2,k\}$ since it is Schreier. This means $\pi_k(X) = A_2$ and $\min \pi_k(X) \geq k+1$, which implies $\pi_k(X) \in \cH(k+1, n)$.
        \item If $k \in A_2$, then $A_1 = \{2\}$ and $k = \min A_2$ by definition. We also note $\pi_k(X) = A_2 \setminus \{k\}$ and because $\abs{X} \geq 3$, we know $n \in \pi_k(X)$, so
        \[
        \min \pi_k(X) \ \geq \ k+1 \ > \ \min A_2 \ \geq \ \abs{A_2} \ > \ \abs{\pi_k(X)},
        \]
        and $\pi_k(X) \in \cH(k+1, n)$.
    \end{enumerate}
    As such, $\pi_k$ is well-defined. We now consider the map
    \[
    \sigma_k: \cH(k+1, n) \ \longrightarrow \ \cB_2^k
    \]
    defined by
    \begin{equation}
        \sigma_k(B) \ = \ B \cup \{2, k\}.
    \end{equation}
    By definition, $n \in B$, so $\abs{\sigma_k(B)} \geq 3$. Note that $\min B \geq k+1$, so
    \[
    \min \sigma_k(B) = \min \{2, k, \min B\} = 2 \quad \text{and} \quad \min (\sigma_k(B) \setminus \{2\}) = k.
    \]
    Since $\sigma_k(B)$ is the union of the two Schreier sets $\{2,k\}$ and $B$, and
    \[
    2 \ \leq \ \min \sigma_k(B),
    \]
    we know $\sigma_k(B) \in \cB_2^k$ by Definition~\ref{def: 2nd order Schreier, nth order Schreier}, so $\sigma_k$ is well-defined.
    
    As it is clear that $\pi_k$ and $\sigma_k$ are inverses of each other, counting $\cB_2^k$ is equivalent to counting $\cH(k+1, n)$, which equals $H(k+1,n)$ by definition. We conclude that
    \begin{align}
        \abs{\cS_2(2,n)} \
        &= \ \abs{\cB_1} + \abs{\cB_2} \nonumber \\
        &= \ \abs{ \bigl\{ \{2,n\} \bigr\}} + \Biggl| {\bigsqcup_{k=3}^{n-1} \cB_2^k} \Biggr| \nonumber\\
        &= \ 1 + \sum_{k=3}^{n-1} \bigl|{\cB_2^k} \bigr| \nonumber \\
        &= \ 1 + \sum_{k=3}^{n-1} H(k+1,n) \nonumber \\
        &= \ 1 + \sum_{k=4}^{n} H(k,n).
    \end{align}
\end{proof}
\begin{remark}
    One consequence of Lemma~\ref{lem: count of min 2 second-order Schreier using H} is that we can easily compute the number of sets in $\cS_2(2,n)$ for $n=4, \ldots, 7$. Specifically, we have
    \[
    \abs{\cS_2(2,n)} \ = \ 1 + \sum_{k=4}^n H(k, n) \ = \ 1 + \sum_{k=4}^n 2^{n-k} \ = \ 2^{n-3},
    \] 
    where the second equality follows from Corollary~\ref{cor: counting min geq a, max b schreier for b<2a}.
\end{remark}
\vspace{1em}
\subsection{Another recurrence relation for $H(a,b)$}~\label{subsec: another recurrence relation for H}
\newline
We now study a recurrence relation for summations of $H(a,b)$. We see that this recurrence relation is consistent with the recurrence relation for $T_n$, as given in Proposition~\ref{prop: recurrence relation for T_n}.
\begin{lemma}\label{lem: recurrence relation for summations of H}
    Let $n \geq 0$. Then, we have
    \begin{equation} \label{eq: recurrence relation for summations of H}
        \sum_{k=4}^{n+6} H(k,n+6) \ = F_{n+4} + \ \sum_{k=4}^{n+5} H(k,n+5) + \sum_{k=4}^{n+4} H(k, n+4)
    \end{equation}
\end{lemma}
\begin{proof}
    Recall from Theorem~\ref{thm: recurrence relation for H} that
    \begin{equation}
    H(a+1, b+2) \ = \ H(a+1, b+1) + H(a, b)
    \end{equation}
    for $a, b \in \bZ_{>0}$ with $a \leq b$. Then,
    \begin{align}
        \sum_{k=4}^{n+5} H&(k, n+5) + \sum_{k=4}^{n+4} H(k, n+4) \nonumber \\
        &= \ H(4, n+5) + \sum_{k=4}^{n+4} \bigl(H(k+1, n+5) + H(k, n+4) \bigr) \nonumber \\
        &= \ H(4, n+5) + \sum_{k=4}^{n+4} H(k+1, n+6) \nonumber \\
        &= \ H(4, n+5) + \sum_{k=5}^{n+5} H(k, n+6).
    \end{align}
    By Theorem~\ref{thm: recurrence relation for H}, we know
    \begin{equation}
        H(4, n+6) \ = \ H(4, n+5) + H(3, n+4),
    \end{equation}
    so to show~\eqref{eq: recurrence relation for summations of H}, it suffices to show
    \begin{equation}\label{eq: reduction of recurrence relation for sums of H}
        H(3,n+4) + 1 \ = \ F_{n+4}.
    \end{equation}
    From Lemma~\ref{lem: counting min geq a, max b Schreier} and~\eqref{eq: counting min geq a, max b schreier, size at least a}, we have
    \begin{equation}
        H(3,n+4) \ = \ \sum_{j=0}^{\min(n+1, 1)} \binom{n+1}{j} + \sum_{\ell = 3}^{\flr{(n+5)/2}} \binom{n+4-\ell}{\ell-1}.
    \end{equation}
    As $n \geq 0$, we know $\min(n+1, 1) = 1$, so
    \begin{equation}
        \sum_{j=0}^{\min(n+1, 1)} \binom{n+1}{j} \ = \ \binom{n+1}{0} + \binom{n+1}{1} \ = \ n+2.
    \end{equation}
    Since
    \[
    1 \ = \ \binom{n+3}{0} \quad \text{and} \quad n+2 \ = \ \binom{n+2}{1},
    \]
    we have
    \begin{align}
        H(3,n+4) + 1 \
        &= \ \sum_{\ell = 3}^{\flr{(n+5)/2}} \binom{n+4-\ell}{\ell-1} + \binom{n+2}{1} + \binom{n+3}{0} \nonumber \\
        &= \ \sum_{\ell = 1}^{\flr{(n+5)/2}} \binom{n+4-\ell}{\ell - 1} \nonumber \\
        &= \ \sum_{k=0}^{\flr{(n+3)/2}} \binom{n+3-k}{k} \nonumber \\
        &= \ F_{n+4}, \label{eq: H(3,n+4) + 1 = F(n+4)}
    \end{align}
    as desired.
\end{proof}
\vspace{1em}
\subsection{Main result}~\label{subsec: proof of main result}
\newline
Recall the statement of Theorem~\ref{thm: count of min 2 second-order Schreier sets}: for $n \geq 0$, we have
\[
\abs{\cS_2(2,n+4)} \ = \ T_n + 1.
\]
We are now ready to prove this theorem.
\begin{proof}
    From Lemma~\ref{lem: count of min 2 second-order Schreier using H}, we have
    \[
    \abs{\cS_2(2,n+4)} \ = \ 1 + \sum_{k=4}^{n+4} H(k, n+4),
    \]
    so, by Definition~\ref{def: convolution of Fibonacci numbers}, it suffices to show 
    \begin{equation} \label{eq: T_n as sum of H(k,n+4)}
        T_n \ = \ \sum_{i=0}^n F_{i+1} F_{n+2-i} \ = \ \sum_{k=4}^{n+4} H(k, n+4).
    \end{equation}
    Note that when $n=0$,
    \begin{equation}
        T_0 \ = \ F_1 F_2 \ = \ 1 \ = \ H(4,4)
    \end{equation}
    and when $n=1$,
    \begin{equation}
        T_1 \ = \ F_1 F_3 + F_2 F_2 \ = \ 2 + 1 \ = \ H(4,5) + H(5,5).
    \end{equation}
    As we have established the base cases, we know from Proposition~\ref{prop: recurrence relation for T_n} that showing
    \begin{equation} \label{eq: WTS for main thm}
        \sum_{k=4}^{n+6} H(k, n+6) \ = \ \sum_{k=4}^{n+5} H(k, n+5) + \sum_{k=4}^{n+4} H(k, n+4) + F_{n+4}
    \end{equation}
    is sufficient to prove~\eqref{eq: T_n as sum of H(k,n+4)}. Since~\eqref{eq: WTS for main thm} is exactly the statement of Lemma~\ref{lem: recurrence relation for summations of H}, we are done.
\end{proof}
\begin{corollary}
    The number of second-order Schreier sets with minimum element $2$ and maximum element $n+4$ that are not (first-order) Schreier is equal to $T_n$.
\end{corollary}
\begin{proof}
    By Definition~\ref{def: Schreier set}, the only element in $\cS_2(2,n+4)$ that is also Schreier is $\{2,n+4\}$. Then, by Theorem~\ref{thm: count of min 2 second-order Schreier sets},
    \begin{equation}
        \# \{ X \in \cS_2(2,n+4) : X \notin \cS_1\} \ = \ T_n.
    \end{equation}
\end{proof}

\section{Problems for Further Investigation}\label{sec: future problems}
A natural question that arises from our work is whether there exists a bijection between the (min $2$, max $n$) second-order Schreier sets and the sequence $T_n$. By Theorem~\ref{thm: count of min 2 second-order Schreier sets}, we observe for fixed $n$ that
\begin{equation}
    |\cS_2(2,n+4)| \ = \ 1 + \bigl(F_1 F_{n+2} + F_2 F_{n+1} + \cdots + F_{n+1} F_2 \bigr),
\end{equation}
and we know by Theorem~\ref{thm: counting Schreier sets} that $F_n = |\cS_1(n)|$. As such, the above equation can be written as
\begin{equation}
    |\cS_2(2,n+4)| \ = \ 1 + \bigl( |\cS_1(1)| \cdot |\cS_1(n+2)| + \cdots + |\cS_1(n+1)| \cdot |\cS_1(2)| \bigr).
\end{equation}
When considering that an arbitrary set in $\cS_2(2,n+4)$ is a union of two Schreier sets, we question whether there exists a bijective map
\[
\Phi: \cS_2(2, n+4) \ \setminus \bigl\{ \{2, n+4 \} \bigr\} \ \longrightarrow \ \bigsqcup_{i=0}^n \cS_1(i+1) \times \cS_1(n+2-i)
\]
involving replacing, shifting, or appending elements that takes unions of Schreier sets to a second-order Schreier set, analogous to the bijective map given in the proof of Theorem~\ref{thm: counting Schreier sets} by Bird~\cite{out-of-the-norm}.

Another direction to explore is whether the count of $\cS_2(a, n)$ for $a > 2$ can be similarly expressed in terms of a known sequence (likely related to the Fibonacci numbers). Since $\cS_2(2,n)$ is an index shift of a convolution of Fibonacci numbers, we question whether $\cS_2(3,n)$ is some manipulation of a double convolution of Fibonacci numbers. More generally, we can also ask whether there exists a formula for counting the $m^{\text{th}}$ Schreier family for $m \geq 2$.

\end{document}